# WEAK MIXING SUSPENSION FLOWS OVER SHIFTS OF FINITE TYPE ARE UNIVERSAL

ANTHONY QUAS AND TERRY SOO

Abstract. Let $S$ be an ergodic measure-preserving automorphism on a non-atomic probability space, and let $T$ be the time-one map of a topologically weak mixing suspension flow over an irreducible subshift of finite type under a Hölder ceiling function. We show that if the measure-theoretic entropy of $S$ is strictly less than the topological entropy of $T$, then there exists an embedding from the measure-preserving automorphism into the suspension flow. As a corollary of this result and the symbolic dynamics for geodesic flows on compact surfaces of negative curvature developed by Bowen [5] and Ratner [31], we also obtain an embedding from the measure-preserving automorphism into a geodesic flow, whenever the measure-theoretic entropy of $S$ is strictly less than the topological entropy of the time-one map of the geodesic flow.

## 1. Introduction

Let $F$ be a homeomorphism of a compact metric space $\Gamma$. Let $S$ be an ergodic measure-preserving automorphism on a non-atomic probability space $\Omega$ with an invariant measure $\nu$. An ***embedding*** of $(\Omega, \nu, S)$ into $(\Gamma, F)$ is a measurable mapping $\Psi : \Omega \to \Gamma$ such that the restriction of $\Psi$ to a set of full measure $\Omega'$ is an injection, and $\Psi(S(\omega)) = F(\Psi(\omega))$ for all $\omega \in \Omega'$. Thus existence of an embedding of $(\Omega, \nu, S)$ into $(\Gamma, F)$ is equivalent to the existence of a $F$-invariant measure $\nu'$ on $\Gamma$ that makes $(\Gamma, \nu', F)$ isomorphic to $(\Omega, \nu, S)$. We say that the topological dynamical system $(\Gamma, F)$ is ***universal*** if for every invertible ergodic measure-preserving system $(\Omega, \nu, S)$ with measure-theoretic entropy strictly less than the topological entropy of $(\Gamma, F)$ there exists an embedding of $(\Omega, \nu, S)$ into $(\Gamma, F)$, and we say that $(\Gamma, F)$ is ***fully universal*** if the embedding can be chosen so that the push-forward of the measure on $\Omega$ is fully supported on $\Gamma$. Let us remark that Lind and Thouvenot [27] used the term 'universal' to mean what we are calling 'fully universal.'

2010 *Mathematics Subject Classification.* Primary 37A35, Secondary 37D40.

Funded in part by NSERC (both authors).

The Krieger finite generator theorem [17, 19] says that the full-shift with a finite number of symbols is universal. Krieger also proved that mixing subshifts of finite type are universal [21] (see also the proof given by Denker [9, Theorem 28.1]) and raised the question of when a homeomorphism of a compact metric is universal. The ***time-one*** map of a flow $(G^t)_{t\in\mathbb{R}}$ is the map $G^1$.

**Theorem 1.** *The time-one map map flow of a geodesic flow on a compact surface of variable negative curvature is universal.*

We will prove Theorem 1 by exploiting the symbolic representation of geodesic flows developed by Bowen [5] and Ratner [31], and proving our main result, that topologically weak mixing suspension flows over irreducible subshifts of finite type under Hölder continuous ceiling functions are also universal.

A related approach was taken by Lind and Thouvenot [27, Section 5] in their proof that hyperbolic (two-dimensional) toral automorphisms are fully universal. They used the fact that hyperbolic toral automorphisms can be represented as irreducible shifts of finite type.

The associated matrix of a quasi-hyperbolic toral automorphism has no roots of unity as eigenvalues, but may have other eigenvalues on the unit circle. Lind and Thouenot [27] asked whether quasi-hyperbolic toral automorphisms (which do not have representations as shifts of finite type [23, Section 6], [24, Theorem 4]) are fully universal. In a forthcoming paper [30], we show that all toral automorphisms are universal and that all quasi-hyperbolic automorphisms are fully universal.

In what follows, we review all the necessary terminology required for stating our main result.

Let $V$ be a finite set of symbols of cardinality $\#V$; we will always assume that $\#V \geq 2$. Let $\theta : V^{\mathbb{Z}} \to V^{\mathbb{Z}}$ be the shift defined by $\theta(y)_i = y_{i+1}$ for all $y \in V^{\mathbb{Z}}$ and all $i \in \mathbb{Z}$. We endow $V^{\mathbb{Z}}$ with the standard product metric and Borel $\sigma$-algebra. Sometimes $V^{\mathbb{Z}}$ is called a full-shift. Let $\mathfrak{A}$ be a square zero-one matrix of order $\#V$. Also assume that $\mathfrak{A}$ is irreducible; that is, for all $1 \leq i, j \leq \#V$, there exists $n \in \mathbb{Z}^+$ such that $\mathfrak{A}^n(i,j) > 0$. Define

$$Y := \left\{ y \in V^{\mathbb{Z}} : \mathfrak{A}(y_i y_{i+1}) = 1 \text{ for all } i \in \mathbb{Z} \right\}.$$

Thus $Y$ is the set of all bi-infinite paths of the directed graph on $V$ with adjacency matrix $\mathfrak{A}$; since $\mathfrak{A}$ is assumed to be irreducible, the graph is strongly connected. We say that $Y$ is an ***irreducible subshift of finite type***; see [25, Chapter 2] for background. More generally, we say that $X$ is a ***subshift***, if $X$ is a closed shift-invariant subset of $V^{\mathbb{Z}}$ (of course $X$ is endowed with the Borel $\sigma$-algebra). We will call a subshift ***non-trivial*** if it does not consist of a finite set of points. (When talking about subshifts, the map will always be the shift map, $\theta$, and it may be left implicit).

Let $Y$ be an irreducible subshift of finite type and let $f : Y \to (0, \infty)$ be Hölder continuous. Set

$$\mathrm{susp}(Y, f) := \{(y, s) \in Y \times [0, \infty) : y \in Y, 0 \leq s < f(y)\}.$$

We define the flow $(T^t)_{t\in\mathbb{R}}$ on $\mathrm{susp}(Y, f)$ by setting $T^t(y, s) = (y, s + t)$ and identifying the points $(\theta y, 0) = (y, f(y))$. More precisely, for all $y \in Y$, set $f^{(0)}(y) := 0$ and $f^{(n+1)}(y) := f^{(n)}(y) + f(\theta^n y)$ for all $n \in \mathbb{Z}$. For all $(y, s) \in \mathrm{susp}(Y, f)$, and $t \in \mathbb{R}$, there is a unique $n \in \mathbb{Z}$ such that $f^{(n)}(y) \leq s + t < f^{(n+1)}(y)$; we set

$$T^t(y, s) = \big(\theta^n y, s + t - f^{(n)}(y)\big).$$

Sometimes $(\mathrm{susp}(Y, f), (T^t)_{t\in\mathbb{R}})$ is called the suspension flow over $Y$ under the ceiling function $f$ or a hyperbolic symbolic flow [4]. The Bowen-Walters distance [8] makes $\mathrm{susp}(Y, f)$ a compact metric space, where a neighbourhood about a point $(y, s) \in \mathrm{susp}(Y, f)$ contains all the points $T^t(w, s) \in \mathrm{susp}(Y, f)$, where $|t|$ is small and $w$ is close to $y$. With respect to the topology generated, $T^t$ is a homeomorphism on $\mathrm{susp}(Y, f)$ for all $t \in \mathbb{R}$.

Following Pollicott and Parry [29, Page 95] we say that the flow $(\mathrm{susp}(Y, f), (T^t)_{t\in\mathbb{R}})$ is ***topologically weak mixing*** if there does not exist a non-constant continuous eigenfunction function $\mathsf{F}$ from the suspension space to the complex unit circle, so that for all $t \in \mathbb{R}$, we have $\mathsf{F} \circ T^t = e^{2\pi \mathbf{i}\beta t}\mathsf{F}$, for some eigenvalue $\beta > 0$; note that for our purposes, the eigenvalue condition is equivalent to the usual definition of topological weak mixing [16, Theorem 2.3 and Theorem 2.5]. (In more general settings, the eigenvalue condition is a weaker condition than the standard condition [16, Theorem 2.3 and Example 2.7]). A simple, yet representative, example in which $(\mathrm{susp}(Y, f), (T^t)_{t\in\mathbb{R}})$ is topologically weak mixing is the Totoki flow given by $Y = \{0, 1\}^{\mathbb{Z}}$, and $f(y) = a_{y_0}$ for all $y \in Y$, where $a_0, a_1$ are positive real numbers whose ratio is irrational [34].

**Theorem 2.** *The time-one map of a topologically weak mixing suspension flow over an irreducible subshift of finite type under a Hölder continuous function is universal.*

We will prove Theorem 2 by first proving a result in the special case of a non-trivial ergodic subshift; Theorem 3 below, in combination with the Jewett-Krieger theorem will imply Theorem 2. See [2, 18, 20, 13] for more on the Jewett-Krieger theorem.

**Theorem 3.** *Let* $(X, \theta)$ *be a non-trivial ergodic subshift with invariant measure* $\mu$. *Let* $(\mathrm{susp}(Y, f), (T^t)_{t\in\mathbb{R}})$ *be a topologically weak mixing suspension flow over an irreducible subshift of finite type* $Y$ *under a Hölder continuous function* $f : Y \to (0, \infty)$. *If the topological entropy of* $(X, \theta)$ *is strictly less than the topological entropy of* $(\mathrm{susp}(Y, f), T^1)$, *then there exists an embedding of* $(X, \mu, \theta)$ *into* $(\mathrm{susp}(Y, f), T^1)$.

Note in Theorem 3, the condition that $(X, \theta)$ has strictly less topological entropy than the topological entropy of $(\mathrm{susp}(Y, f), T^1)$ implies that the measure-theoretic entropy of $(X, \mu, \theta)$ is strictly less than the topological entropy of $(\mathrm{susp}(Y, f), T^1)$.

*Proof of Theorem 2.* By the Jewett-Krieger theorem [20], a non-atomic probability space endowed with an ergodic measure-preserving automorphism with finite measure-theoretic entropy is isomorphic to a uniquely ergodic subshift equipped with its invariant measure. Notice that by the variational principle [12], the topological entropy of the uniquely ergodic subshift is equal to the measure-theoretic entropy of the initial automorphism. Hence, composing this isomorphism with the embedding given by Theorem 3 produces the required embedding. □

Let us remark that the embedding that we define to prove Theorem 3 is not continuous.

**Question 1.** *Let* $(X, \theta)$ *be a subshift equipped with an invariant measure* $\mu$. *Let* $(\mathrm{susp}(Y, f), (T^t)_{t\in\mathbb{R}})$ *be a topologically weak mixing suspension over an irreducible subshift of finite type* $Y$ *under a Hölder continuous function* $f : Y \to (0, \infty)$.

*If the measure-theoretic entropy of* $(X, \theta)$ *is (strictly) less than the topological entropy of* $(\mathrm{susp}(Y, f), T^1)$, *under which conditions must there exist an embedding* $\Psi : X \to \mathrm{susp}(Y, f)$ *where the continuity points have full measure?*

Mappings satisfying the property in the question are said to be *finitary* (see [32] and [15] for more information).

It remains to prove Theorems 1 and 3, but next we discuss an application of Theorem 1 that was the original motivation of Theorem 2. We thank Jean-Paul Thouvenot [33] who suggested that Theorem 2 could be used to prove Theorem 1, from which one can obtain a negative answer to the following problem posed by François Ledrappier, Federico Rodriguez Hertz, and Jana Rodriguez Hertz [22]:

> Let $(G^t)_{t\in\mathbb{R}}$ be geodesic flow on a compact surface of variable negative curvature $\mathcal{M}$, with unit tangent bundle $\mathrm{UT}(\mathcal{M})$. Suppose that $G^n(w) \neq w$ for all $n \in \mathbb{Z}\backslash\{0\}$

> and for all $w \in \mathrm{UT}(\mathcal{M})$. If a measure $\lambda$ is invariant and ergodic under the group $(G^n)_{n\in\mathbb{Z}}$, must it be invariant under $(G^t)_{t\in\mathbb{R}}$?

**Corollary 4.** *Let $(G^t)_{t\in\mathbb{R}}$ be geodesic flow on a compact surface of variable negative curvature. There is a measure $\lambda$ that is invariant and ergodic under the group $(G^n)_{n\in\mathbb{Z}}$, but not invariant under $(G^t)_{t\in\mathbb{R}}$.*

Corollary 4 is immediate from Theorem 1, and the simple fact that there exists an ergodic measure-preserving automorphism that has arbitrarily small entropy and does not admit a square root. Let us remark that by appealing to a result of Ornstein [28], which says that there exists a $K$-automorphism that does not admit a square root, one can also require that $(UT(\mathcal{M}), \lambda, T^1)$ is a $K$-automorphism.

The following related question is also due to François Ledrappier, Federico Rodriguez Hertz, and Jana Rodriguez Hertz [22]:

> Let $(G^t)_{t\in\mathbb{R}}$ be geodesic flow on a compact surface of variable negative curvature $\mathcal{M}$, with unit tangent bundle $\mathrm{UT}(\mathcal{M})$. Suppose that $G^n(w) \neq w$ for all $n \in \mathbb{Z}\backslash\{0\}$ and for all $w \in \mathrm{UT}(\mathcal{M})$. If a set is minimal under the action of $(G^t)_{t\in\mathbb{R}}$, must it be minimal under $(G^n)_{n\in\mathbb{Z}}$?

In the appendix, we give a negative answer to this question (we note that this implies a negative answer to the earlier question).

We prove Theorem 1 and Corollary 4 in Section 5. In the next section, we give an outline of the proof of Theorem 3 that will also give criteria for topological weak mixing for suspensions flows of subshifts of finite type under Hölder continuous ceiling functions. In Section 3, we will assemble some lemmas that will help us define the embedding of Theorem 3 in Section 4.

## 2. Background and proof sketch

In this section, we first introduce some basic terminology. Second, we will discuss topological weak mixing for suspension flows and thirdly, we will discuss the basic approach and highlight the main ideas of the proof.

2.1. **Basic terminology.** Let $X$ be a non-trivial subshift. Let $x \in X$ and $n \in \mathbb{Z}^+$. We say that the finite string $x_0 \cdots x_{n-1} = x_0^{n-1}$ is a ***block of size* n**. Let $B^n(X) := \left\{x_0^{n-1} : x \in X\right\}$ denote the set of all blocks of size $n$ in $X$. We refer to all members of $B(X) = \bigcup_{n\in\mathbb{Z}^+} B^n(X)$, as **X-*blocks***. Given two $X$-blocks $x_0^{n-1}$ and $z_0^{m-1}$ we let their concatenation be given by

$$x_0 \dots x_{n-1} z_0 \dots z_{m-1} = x_0^{n-1} z_0^{m-1}, \tag{1}$$

whenever the concatenation is also a $X$-block. We say that a block $x_0^{n-1}$ ***appears*** in $z_0^{m-1}$ or $z \in X$ if there exists $k \in \mathbb{Z}$ such that $z_k^{k+n-1} = x_0^{n-1}$. In order to distinguish between $x \in X$ and $\theta x \in X$, sometimes we will use the symbol "**.**" to indicate the position of the origin, so that if $x = \cdots x_{-2} x_{-1} . x_0 x_1 x_2 \cdots$ then $\theta x = \cdots x_{-2} x_{-1} x_0 . x_1 x_2 \cdots$.

Let $Y$ be an irreducible subshift of finite type, and let $f : Y \to (0, \infty)$. Given $(y, s) \in \mathrm{susp}(Y, f)$ for each $m \in \mathbb{Z}$, there exists a $t \in \mathbb{R}$ such that $T^t(y, s) = (\theta^m y, 0)$; for any $n > m$, we say that the $Y$-block $y_m^{n-1}$ ***begins at t (in (y, s))***. For $y \in Y$, we define the ***length in the suspension*** of the $Y$-block $y_m^{n-1}$ in $y$ by

$$\mathrm{Len}(y, f; m, n-1) = f^{(n-m)}(\theta^m y) = \sum_{i=m}^{n-1} f(\theta^i y) \tag{2}$$

Note the difference between the size of a block (number of symbols) and length of a block in the suspension. In general, the length in the suspension of any $Y$-block $y_m^n$ may depend on the whole bi-infinite sequence $y$. We define the ***maximum length in the suspension*** of a block $y_0^{m-1}$ by

$$\mathrm{MLen}(y_0^{m-1}, f) := \max \left\{ f^{(m)}(z) : z \in Y, z_0^{m-1} = y_0^{m-1} \right\}.$$

Similarly, for a block $y_0^{m-1}$ define

$$\mathrm{mLen}(y_0^{m-1}, f) := \min \left\{ f^{(m)}(z) : z \in Y, z_0^{m-1} = y_0^{m-1} \right\}.$$

Let $t > 0$, for $Y' \subseteq Y$ define

$$B_t(Y', f) := \{ B \in B(Y') : \mathrm{MLen}(B, f) \le t \}. \tag{3}$$

We say that $(\mathrm{susp}(Y, f), (T^t)_{t \in \mathbb{R}})$ satisfies ***specification*** if for all $\varepsilon > 0$, there exists $L_\varepsilon$ such that given any two blocks $A_0, A_1 \in B(Y)$, and a real number $L > L_\varepsilon + \mathrm{MLen}(A_0, f)$, there exists $(y, 0) \in \mathrm{susp}(Y, f)$ such that $A_0$ begins at 0 and $A_1$ begins within $\varepsilon$ of $L$. Consider again the Totoki flow where $Y = \{0, 1\}^{\mathbb{Z}}$, and $f(y) = a_{y_0}$ for all $y \in Y$, where $a_0, a_1$ are positive real numbers. It is easy to verify that $(\mathrm{susp}(Y, f), (T^t)_{t \in \mathbb{R}})$ satisfies specification if and only if $a_0/a_1$ is irrational.

2.2. **Topological weak mixing.** We discuss in this sub-section some criteria equivalent to topological weak mixing of suspension flows.

Let $Y$ be a subshift of finite type, and let $f : Y \to (0, \infty)$. If $z \in Y$ is periodic, we let $\mathrm{per}(z)$ be the least positive integer such that $\theta^{\mathrm{per}(z)} z = z$. The ***period in the suspension*** is

$$\mathrm{PLen}(z, f) := \mathrm{Len}(z, f; 1, \mathrm{per}(z)).$$

For any non-empty $A \subset \mathbb{R} \setminus \{0\}$, let $\operatorname{span}(A)$ denote the closure of the finite integer combinations of elements of $A$. Since this is a closed subgroup of $\mathbb{R}$, it is either $\mathbb{R}$ or a subgroup of the form $c\mathbb{Z}$ for some $c > 0$. In the latter case, we define $\gcd(A) = c$ and if $\operatorname{span}(A) = \mathbb{R}$, we define $\gcd(A) = 0$.

We say that $f : Y \to (0, \infty)$ is ***filling*** if it satisfies the following properties.

(i) The function $f$ is Hölder continuous.
(ii) The set of all periods in the suspension has greatest common divisor equal to zero.

We say that two functions $f$ and $g$ are ***Hölder cohomologous*** if there exists a Hölder continuous function $h$ such that $f = g + h - h \circ \theta$.

**Proposition 5** (Equivalent notions of topological weak mixing)**.** *Let* $(\operatorname{susp}(Y, f), (T^t)_{t\in\mathbb{R}})$ *be the suspension flow over an irreducible subshift of finite type* $Y$ *under a Hölder continuous ceiling function* $f : Y \to (0, \infty)$*. The following conditions are equivalent.*

*(a) The suspension flow is topologically weak mixing.*
*(b) The ceiling function is filling.*
*(c) The suspension flow satisfies specification.*
*(d) The ceiling function* $f$ *is not Hölder cohomologous to a function* $h$ *taking values in* $\beta\mathbb{Z}$ *for some* $\beta > 0$*.*

This proposition is essentially standard. We remark that our proofs of Theorems 3 and Theorem 1 only require Proposition 5 ((a) $\Rightarrow$ (b)); we give brief details of this implication. A slightly generalized version of (b)$\Rightarrow$ (c) appears in our proof of Lemma 8(F).

To see (a) implies (b), we argue by the contrapositive. If $f$ fails to be filling, the periods in the suspension must all lie in a subgroup $\beta\mathbb{Z}$ of the reals. In this case a version of the Livschitz theorem appearing in the book of Parry and Pollicott [29, Proposition 5.2] implies that $f$ is cohomologous to a function $h$ taking values in $\beta\mathbb{Z}$ (so that (d) fails). This can then be used to construct a continuous eigenfunction

$$\mathsf{F}(y, s) = e^{2\pi \mathbf{i}(s-h(y))/\beta}$$

which satisfies $\mathsf{F}(T_t(y, s)) = e^{2\pi \mathbf{i} t/\beta}\mathsf{F}(y, s)$, contradicting (a).

Proposition 5 explains the role of the assumption that the suspension flow is topologically weak mixing in Theorem 3. If the suspension flow is not assumed to be topologically weak-mixing, the conclusion of Theorem 3 may fail. If Condition (d) fails, then the ceiling function is cohomologous to a function that takes values in $\beta\mathbb{Z}$. If $\beta$ is not of the form $1/n$, then any measure-preserving transformation that can be embedded into $\text{susp}(Y, f)$ must have a (measurable) eigenfunction and hence cannot be weak-mixing (in the measure-theoretic sense). On the other hand if $\beta$ is of the form $1/n$, then more straightforward methods establish the existence of an embedding.

2.3. **The basic approach.** Suppose that $(X, \theta)$ is a non-trivial ergodic subshift with invariant measure $\mu$. Let $(\text{susp}(Y, f), (T^t)_{t\in\mathbb{R}})$ be a topologically weak mixing suspension over an irreducible subshift of finite type $Y$ under a Hölder continuous function $f : Y \to (0, \infty)$. Assume that the topological entropy of $(X, \theta)$ is strictly less than topological entropy of $(\text{susp}(Y, f), T^1)$.

Our approach is to partition a point $x \in X$ into $X$-blocks $\left(x_{n_i}^{n_{i+1}-1}\right)_{i\in\mathbb{Z}}$, and then encode each of these blocks into corresponding $Y$-blocks $\left(\phi\left(x_{n_i}^{n_{i+1}-1}\right)\right)_{i\in\mathbb{Z}}$. Specification allows us to piece together $Y$-blocks and produce an element $\Psi(x) = (y(x), s(x)) \in \text{susp}(Y, f)$; however this must be done is such way so as to guarantee that $\Psi(\theta(x)) = T^1(\Psi(x))$. In addition, some care is needed to ensure that we can recover $x$ from $\Psi(x)$.

The following special case of Alpern's multiple Rokhlin tower theorem [1, 10] will be used to produce the partitions.

**Lemma 6** (Alpern's theorem [1]). *Let $(X, \theta)$ be a non-trivial ergodic subshift with an invariant measure $\mu$. For any positive integer $n$ there exist two measurable sets $Q_1, Q_2$, such that $\{\theta^j(Q_1)\}_{j=0}^{n-1} \cup \{\theta^j(Q_2)\}_{j=0}^{n}$ give a partition of $X$ (modulo a null set).*

Using Lemma 6, we can define a function

$$M(x) = \{u \in \mathbb{Z} : \theta^u(x) \in Q_1 \cup Q_2\} \subset \mathbb{Z}$$

that gives an equivariant subset of $\mathbb{Z}$; that is, for $x \in X$, if $M(x)$ is the subset of $\mathbb{Z}$ assigned to $x$, then $M(\theta x)$ is obtained by subtracting 1 from each elements of $M(x)$. Let $(n_i)_{i\in\mathbb{Z}} = (n_i(x))_{i\in\mathbb{Z}}$ be the enumeration of $M(x)$ such that $\ldots < n_{-2}(x) < n_{-1}(x) < n_0(x) \leq 0 < n_1(x) < n_2(x) < \ldots$. The point $x$ can then be partitioned into $X$-blocks as

$$x = \cdots\, x_{n_{-1}}^{n_0-1}\, x_{n_0}^{n_1-1}\, x_{n_1}^{n_2-1} \cdots .$$

As a consequence of the fact that the topological entropy of $(X, \theta)$ is strictly less than that topological entropy of $T^1$, there will exist $\alpha \in (0, 1)$ such that for all $n$ sufficiently large, we have

$$\#(B^n(X) \cup B^{n+1}(X)) < \#B_{(1-\alpha)n}(Y, f);$$

see the proof of Lemma 8 (H) for details. Thus there exists an injection from $\phi_1 : B^n(X) \cup B^{n+1}(X) \to B_{(1-\alpha)n}(Y, f)$, which encodes $n$-blocks and $(n+1)$-blocks of $X$ into (shorter) $Y$-blocks of length in the suspension at most $(1-\alpha)n$.

For each $i \in \mathbb{Z}$, we will set $A_i = x_{n_i}^{n_{i+1}-1}$ and

$$\psi_1(x) = \cdots . \phi_1(A_0) \cdots \in Y,$$

where $\phi_1(A_0)$ is suitably extended to become an element of $Y$. Let $\Psi_1(x) = T^{-n_0(x)}(\psi_1(x), 0)$. At this stage $\Psi_1$ encodes only the block $A_0$, but we will define subsequent $\Psi_j$ that will encode successively more of the blocks $A_i$ that appear in $x$. Observe that since the $n_i(x)$ were chosen in an equivariant way, we have that $\Psi_1(\theta^k x) = T^k \Psi_1(x)$ for all $k \in [n_0, n_1)$. In this way, $\Psi_1$ maps $x$ in such a way that the block $A_0$ is encoded and starts *exactly* at time $n_0(x)$ in the image. The above will ensure that the encoded block has length in the suspension of approximately $(1-\alpha)n$, so that in the suspension, the encoded version of the block is approximately $\alpha n$ shorter than the original block.

By an additional application of Alpern's multiple Rokhlin tower theorem, we may choose a much sparser (two-sided) subsequence $(n_i')_{i\in\mathbb{Z}}$ of $(n_i)$ in an equivariant way, where again $\ldots n_{-1}' < n_0' \leq 0 < n_1' \ldots$.

Define level two blocks by $A_i' = x_{n_i'}^{n_{i+1}'-1}$, so that $A_0'$ consists of the concatenation of $A_a, A_{a+1}, \ldots, A_{b-1}$ for some $a \leq 0 < b$. As long as the size $n$ above (appearing in the first application of Alpern's lemma) is taken sufficiently large (depending on $\alpha$ and how long it takes for specification to 'kick in'), it is possible to place filler blocks in the gaps of length in the suspension approximately $n\alpha$ between the encodings of the $(A_i)_{a\leq i<b}$ to make a second order block in which each encoded block starts within a fixed precision, $2^{-1}$ say, of the desired starting location $n_i(x)$. Again this finite block is extended to a point $\psi_2(x)$ in $Y$. The second approximation to the embedding is then given by $\Psi_2(x) = T^{-n_0'(x)}(\psi_2(x), 0)$. Now the encoded $A_0'$ block starts exactly at $n_0'(x)$, while the $(A_i)_{a\leq i<b}$ start approximately at $n_i(x)$.

In order to iterate this construction, one minor issue is that to achieve the packing at the next stage, it is necessary to ensure again that the encoded level two blocks are shorter (by an amount that depends on the placement accuracy that we want to obtain in the third level) than the level two blocks in the source. We achieve this by using shorter filler between some of the last level one blocks (the *crumple zone*) forming the level two block.

The inductive step is illustrated in Figure 1.

The Borel-Cantelli lemma will be used to ensure that almost every point only is in finitely many crumple zones. For a point that it is in finitely many crumple zones, we will check that $\Psi_n(x)$ forms a Cauchy

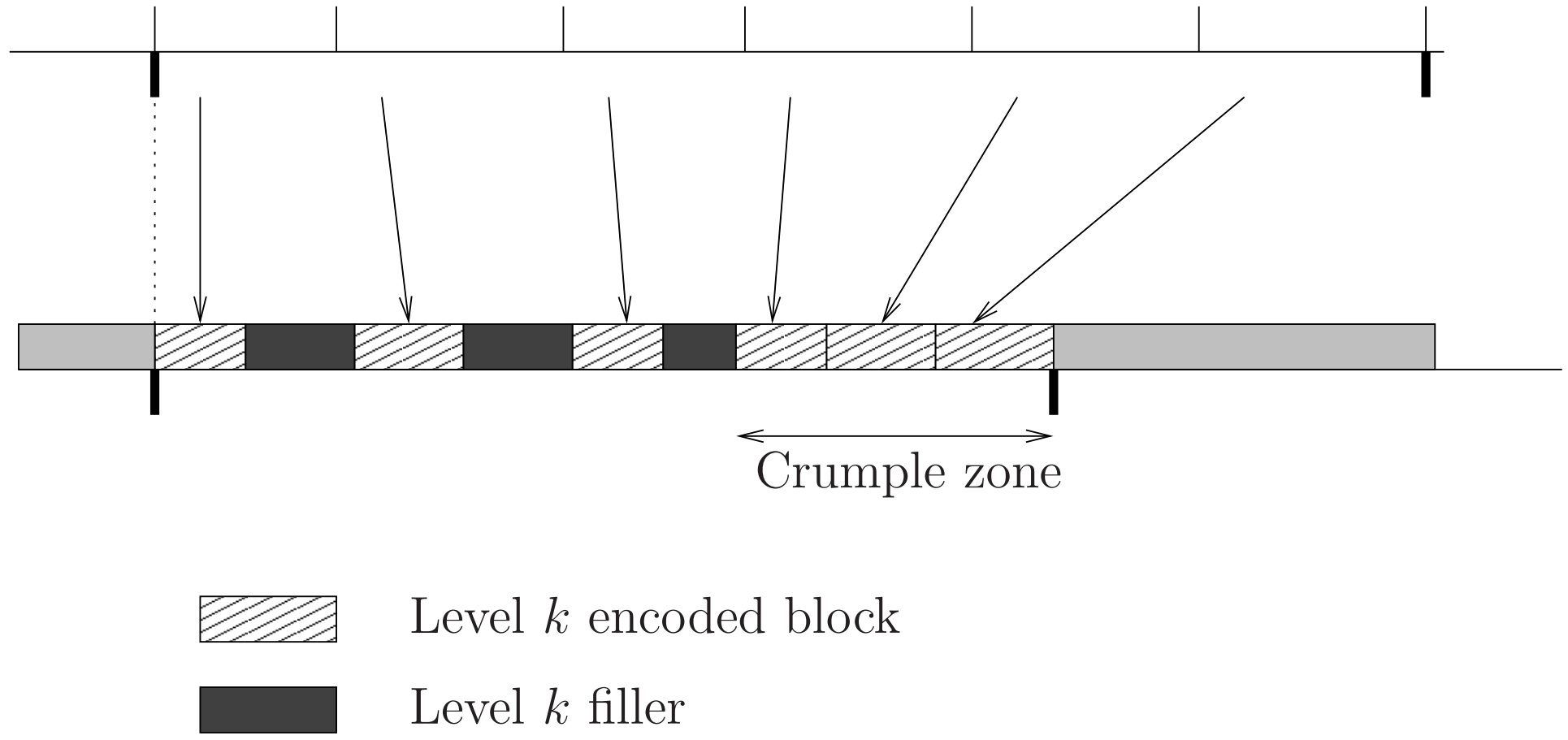

FIGURE 1. Inductive step building the $(k+1)$st level map from the $k$th level map. Notice the crumple zone on the right is compressed to ensure the image block is shorter than the source, allowing subsequent levels to accurately position blocks

sequence, and therefore converges to a limit. This will allow us to define the embedding.

More formally, we will obtain a sequence of mappings from $X$ to $\mathrm{susp}(Y,f)$ with the following properties.

**Proposition 7.** *Suppose that $(X,\theta)$ is a non-trivial ergodic subshift with invariant measure $\mu$. Let $(\mathrm{susp}(Y,f),(T^t)_{t\in\mathbb{R}})$ be a topologically weak mixing suspension flow over an irreducible subshift of finite type $Y$ under a Hölder continuous function $f: Y\to(0,\infty)$. Let $W$ be an arbitrary $Y$-block. If the topological entropy of $(X,\theta)$ is strictly less than the topological entropy of $(\mathrm{susp}(Y,f),T^1)$, then there exists a sequence of measurable maps $\Psi_i: X\to\mathrm{susp}(Y,f)$ and a subset $X'\subset X$ with $\mu(X')=1$, with the following properties.*

*(a) For all $x\in X'$ and for all sufficiently large $i\in\mathbb{Z}^+$ we have*
$$\Psi_i(\theta x)=T^1\Psi_i(x).$$
*(b) For all $x\in X'$, we have $\lim_{i\to\infty}\Psi_i(x):=\Psi(x)=(\mathbf{y}(x),\mathbf{s}(x))\in\mathrm{susp}(Y,f)$ exists.*

*(c) For all $x\in X'$, we have that $\mathbf{y}(x)$ is aperiodic; that is, there does not exist $m\in\mathbb{Z}\setminus\{0\}$ such that $\theta^m\mathbf{y}(x)=\mathbf{y}(x)$.*

*(d) For all $x\in X'$, we have that $\mathbf{y}^{-1}(\mathbf{y}(x))\subseteq\{\theta^n x: n\in\mathbb{Z}\}$.*

*(e) For all $x\in X'$, the block $W$ appears in $\mathbf{y}(x)$ infinitely often in the positive coordinates and in the negative coordinates.*

The proof of Theorem 3 follows easily from Proposition 7. Condition (e) is not required for the proof of Theorem 3, but is necessary for our proof of Theorem 1.

*Proof of Theorem 3.* We apply Proposition 7 with an arbitrary $Y$-block $W$. That $\Psi(\theta(x)) = T^1(\Psi(x))$ follows from (a) and (b). To see that $\Psi$ is an injection, let $x \neq x' \in X'$ and assume $\Psi(x) = \Psi(x')$. By Property (d), $x' = \theta^n x$ for some $n \in \mathbb{Z}$. This implies $\Psi(x) = \Psi(x') = T^n(\Psi(x))$ which gives a contradiction using (c).

□

## 3. Key Ingredients

In this section, we will define, using Alpern's multiple Rokhlin tower theorem [1, 10], equivariant subsets of $\mathbb{Z}$, as a function on a non-trivial ergodic subshift; following Keane and Smorodinsky [14, 15], we think of elements of these subsets as *markers.* These two ingredients will allow us to define the embedding of Theorem 3 in Section 4.

Before we state Lemma 8, we need some additional notation. Let $Y$ be an irreducible subshift of finite type. Given any $A, B \in B(Y)$ there is a $C \in B(Y)$ such that $ACB \in B(Y)$. In the case, where $AB \notin B(Y)$, we define the ***connecting block*** of $A$ and $B$ to be the block $C[A, B]$ such that $A \cdot C[A, B] \cdot B \in B(Y)$ with the least maximum length in the suspension, where we break ties by assigning a lexicographic order on $B(Y)$; in the case where $AB \in B(Y)$, we let $A\, C[A, B]\, B = AB$, and say that $C[A, B]$ is an empty connecting block, and assign it zero maximum length in the suspension. Sometimes to make expressions more readable, we will omit a connecting block and write

$$A\ C[A, B]\ B = A \frown B.$$

Clearly, the maximum lengths in the suspension of all connecting blocks are bounded above.

Given a block $A \in B(Y)$, we let

$$Y \setminus A := \{y \in Y :\ A \text{ does not appear in } y\}.$$

Note that $Y \setminus A$ is a subshift of finite type.

For any $B \in B(Y)$ such that $BB \in B(Y)$, let $B^k := \overbrace{B \cdots B}^{k \text{ times}}$ for any integer $k \geq 1$, and let $\overline{B} \in Y$ be the periodic point with

$$\overline{B} = \cdots BB.BB \cdots .$$

**Lemma 8.** *Let* $(\mathrm{susp}(Y, f), (T^t)_{t\in\mathbb{R}})$ *be a topologically weak mixing suspension flow over an irreducible subshift of finite type* $Y$ *under a Hölder* *continuous ceiling function $f$. Let $(X,\theta)$ be a non-trivial subshift. Suppose that the topological entropy of $(X,\theta)$ is strictly less than the topological entropy of $(\mathrm{susp}(Y,f),T^1)$.*

*Given any $W \in B(Y)$, there exist $Y$-blocks $\langle, \rangle, \mathbf{L}, \mathbf{R}_1, \mathbf{R}_2 \ldots$, and non-decreasing integers $(K_n)_{n\in\mathbb{Z}^+}, (k_n)_{n\in\mathbb{Z}^+}, (L_n)_{n\in\mathbb{Z}^+}$ with the following properties.*

*(A)* *The concatenations $\mathbf{L}\,\mathbf{L}$ and $\mathbf{R}_i\,\mathbf{R}_i$ for all $i \in \mathbb{Z}^+$, are also $Y$-blocks.*

*(B)* *We have $\gcd\big(\mathrm{PLen}(\overline{\mathbf{L}}, f), \mathrm{PLen}(\overline{\mathbf{R_i}}, f)\big) \le 2^{-i}$ for every $i \ge 1$.*

*(C)* *The block $W$ appears in the blocks $\langle$ and $\rangle$.*

*(D)* *The concatenations $\langle\,\langle$ and $\rangle\,\rangle$ are also $Y$-blocks. The block $\langle$ does not appear in the periodic point $\overline{\rangle}$, and the block $\rangle$ does not appear in the periodic point $\overline{\langle}$.*

*(E)* *For all nonnegative integers $\ell$ and $r$, the $Y$-block*

$$C[\langle, \mathbf{L}]\ \mathbf{L}^{\ell} \frown \mathbf{R}_i^r\ C[\mathbf{R}_i, \rangle]$$

*does not contain $\langle^{k_i}$ or $\rangle^{k_i}$.*

*(F)* *For any two blocks $A_0, A_1 \in B(Y)$ and for all $L > L_i + \mathrm{MLen}(A_0, f)$, there exists $(z,0) \in \mathrm{susp}(Y,f)$ with the following properties.*

*(i)*

$$z = \cdots \rangle\,\rangle \frown \mathbf{.}\, A_0\ D_0\ A_1 \frown \langle\,\langle \cdots,$$

*where $D_0$ is given by*

$$C[A_0, \langle]\ \langle^{K_i} \frown \mathbf{L}^{\ell} \frown \mathbf{R}_i^r \frown \rangle^{K_i}\ C[\rangle, A_1]\,, \tag{4}$$

*for some positive integers $\ell$ and $r$.*

*(ii)* *The block $A_1$ begins within $2 \cdot 2^{-i}$ of $L$.*

*(iii)* *For all $w \in Y$, and all positive integers $\ell'$ and $r'$, we have that the block $A_1$ begins within $3 \cdot 2^{-i}$ of $L$ in $(w', 0)$, where*

$$w' = \cdots\cdots \rangle\,\rangle \frown \mathbf{.}\, A_0\ D_0\ A_1\ D\, w_0 w_1 \cdots,$$

*and $D$ is the block*

$$C[A_1, \langle]\ \langle^{K_i} \frown \mathbf{L}^{\ell'} \frown \mathbf{R}_i^{r'} \frown \rangle^{K_i}\ C[\rangle, w_0]$$

*(G)* *The sequence $(K_n)_{n\in\mathbb{Z}}$ satisfies*

$$K_n > 2K_{n-1} + 2k_{n-1} + n. \tag{5}$$

*(H)* *There exists $\alpha \in (0,1)$ such that for all $i$ sufficiently large, there exists $\mathbf{F} \in B(Y)$ and an injection*

$$\phi : B^i(X) \cup B^{i+1}(X) \to B_{(1-\alpha)i}(Y \setminus \mathbf{F}, f)$$

*such that the following hold.*

*(i)* *The block $\mathbf{F}$ appears in both $\langle$ and $\rangle$.*

*(ii)* *For all $B \in B^i(X) \cup B^{i+1}(X)$, we have $\rangle\ \phi(B)\ \langle \in B(Y)$.*

*(iii) For all $B \in B^i(X) \cup B^{i+1}(X)$, the only appearance of $\rangle$ in $\rangle\, \phi(B)\, \langle$ is in the leftmost position; and the only appearance of $\langle$ is in the rightmost position.*

The map from Lemma 8 (H) allows us to encode $X$-blocks into $Y$-blocks, which we think of as blocks of information. The resulting blocks of information are extended to become elements of $Y$ using the blocks $\langle$ and $\rangle$, and then joined together using Lemma 8 (F), which allows us to prescribe where they begin within a certain accuracy, at the cost of placing a filler block in between two blocks of information. Thus we obtain an alternating sequence of information and filler blocks.

We need to be able to distinguish the information from the filler, and this is done by bracketing the filler with concatenations of the blocks $\langle$ and $\rangle$. One technicality is that although we can demand that the blocks $\langle$ and $\rangle$ do not appear in the information, it turns out that it is too much to ask that they do not appear in the filler; see Example 1 in Section 6. However, we can demand that the number of $\langle$ and $\rangle$ that appear in the filler is smaller than the number of $\langle$ and $\rangle$ blocks used to bracket the filler. Although, the blocks $\langle$ and $\rangle$ do not appear in the information, we need condition (Hiii) of Lemma 8, to assure us that we are able to distinguish the beginning or end of a block of information with the end of beginning of a $\rangle$ or $\langle$, respectively.

*Proof of Lemma 8* (A) *and* (B). We first select the blocks $\mathbf{L}$ and $(\mathbf{R}_i)$. We consider two cases. If there exist $\mathbf{L}$ and $\mathbf{R}$ for which $\mathbf{L}\,\mathbf{L}$ and $\mathbf{R}\,\mathbf{R}$ are both legal and such that $\gcd(\mathrm{PLen}(\overline{\mathbf{L}}, f), \mathrm{PLen}(\overline{\mathbf{R}}, f)) = 0$, then let $\mathbf{R}_i = \mathbf{R}$ for each $i$.

Otherwise, let $z^1 \in Y$ be a periodic point. Let $\mathbf{L}$ be the periodic block forming $z^1$. Fix an integer $i \geq 1$. There exists a positive integer $n$ so that $\mathrm{PLen}(z^1, f)/n < 2^{-i}$. By Proposition 5 (a) $\to$ (b), we have that $f$ is filling; thus there exists a periodic point $z^2 \in Y$ such that $\mathrm{PLen}(z^2, f) \notin (\mathrm{PLen}(z^1, f)/n!)\mathbb{Z}$. Let $\mathbf{R}_i$ be the periodic block forming $z^2$. Let $\gcd\big(\mathrm{PLen}(z^1, f), \mathrm{PLen}(z^2, f)\big) = r$. Since $r > 0$, but is not a multiple of $\mathrm{PLen}(z^1, f)/n!$, we must have $\mathrm{PLen}(z^1, f) = rm$ for some integer $m > n$, from which we conclude that $r < 2^{-i}$.

Notice that the periods in the suspension of the $\overline{\mathbf{R}_i}$ must approach infinity as $i \to \infty$ as otherwise infinitely many $\overline{\mathbf{R}_i}$ would be identical by the pigeonhole principle and could be expressed as $\overline{\mathbf{R}}$. This would fall into the case considered above. $\square$

*Proof of Lemma 8* (C) *and* (D). Let $V = W \cdot C[W, W]$ be the shortest word containing $W$ that can be periodically concatenated. Fix a $p \geq 2$ such that $h_{\mathrm{top}}(X, \theta) < h_{\mathrm{top}}(\mathrm{susp}(Y \setminus V^p, f), T^1)$ (such a $p$ exists as the right side converges to $h_{\mathrm{top}}(\mathrm{susp}(Y,f),T^1)$ as $p\to\infty$ [26]). Set $\mathbf{F}=V^p$.

Let $\langle$ and $\rangle$ be two blocks of the form $\mathbf{F}^{10}U\mathbf{F}^{10}$, where the $U$'s are chosen to ensure that $\langle$ does not appear in $\overline{\rangle}$, $\rangle$ does not appear in $\overline{\langle}$ and neither of $\langle$ or $\rangle$ can be written as a power of a shorter block.

To see why such $U$'s must exist, we argue as follows. Choose a $Q'\in B(Y)$ that does not appear in $\overline{\mathbf{F}}$. Choose a $Q\in B(Y)$ that does not appear in

$$\cdots\mathbf{F}\mathbf{F}\frown Q'\frown\mathbf{F}\mathbf{F}\cdots.$$

Set

$$\langle := \mathbf{F}^{10}\frown Q\frown \mathbf{F}^{10}.$$

Choose an integer $n>10$ such that $\mathbf{F}^n$ does not appear in $\overline{\langle}$ and so that $\mathrm{mLen}(\mathbf{F}^n,f)>\mathrm{MLen}(Q,f)$. Set

$$\rangle := \mathbf{F}^{n}\frown Q'\frown \mathbf{F}^{n}.$$

Finally, since $\mathbf{F}^n$ does not appear in $\overline{\langle}$, it follows that $\rangle$ does not appear in $\overline{\langle}$. By the definition of $Q$ and the fact that $\mathrm{MLen}(Q,f)<\mathrm{mLen}(\mathbf{F}^n,f)$, we have that $Q$ does not appear in $\overline{\rangle}$, and thus neither does $\langle$.

We may additionally assume that $p$ was chosen so that

$$\mathrm{PLen}(\overline{\langle},f)>\mathrm{mLen}(\mathbf{F},f)=\mathrm{mLen}(V^p,f)>\mathrm{PLen}(\overline{\mathbf{L}},f)$$

and also $\mathrm{PLen}(\overline{\rangle},f)>\mathrm{PLen}(\overline{\mathbf{R}},f)$ in the case that $\mathbf{R}_i=\mathbf{R}$ for all $i\geq 1$. □

*Proof of Lemma 8* (E). Observe that (E) is satisfied for some $(k_i)$ if and only if $\overline{\mathbf{L}}$ and $\overline{\mathbf{R}_i}$ are not in the orbit of $\overline{\langle}$ or $\overline{\rangle}$ for all $i\in\mathbb{Z}^+$. For this it suffices to check that the periods in the suspension of $\overline{\mathbf{L}}$ and $\overline{\mathbf{R}_i}$ differ from those of $\overline{\langle}$ and $\overline{\rangle}$.

Since $\overline{\langle}$ and $\overline{\rangle}$ have longer period in the suspension than $\overline{\mathbf{L}}$ by construction, the conclusion is straightforward for $\overline{\mathbf{L}}$. Similarly for the $\overline{\mathbf{R_i}}$ if all of the $\mathbf{R}_i$ are equal to $\mathbf{R}$.

Otherwise, by the remark made in the proof of parts (A) and (B), for sufficiently large $i$, the period in the suspension of $\overline{\mathbf{R}_i}$ would exceed the periods in the suspension of $\overline{\langle}$ and $\overline{\rangle}$. Discarding those $\mathbf{R}_i$ for which this fails and re-numbering (this can be done without affecting any previous stage of the construction), we obtain the desired conclusion. □

*Proof of Lemma 8* (F) *and* (G). It will become clear that it does not matter how we define $K_i$ for the purposes of conditions (Fi) and (Fii). We will need to make $K_i$ large to satisfy condition (Fiii). Hence satisfying (5) is trivial.

A straightforward calculation shows that for all $v, v' \in Y$, and $-\iota(v,v') \le i \le j \le \iota(v,v')$, we have that

$$|\operatorname{Len}(v, f; i, j) - \operatorname{Len}(v', f; i, j)| \le M2^{-c[\iota(v,v')-\max(|i|,|j|)]}, \qquad (6)$$

where $\iota(v,v') = \min\left\{|j| \in \mathbb{Z} : v_j \ne v'_j\right\}$, $M \ge \|f\|$ is some constant independent of $v, v'$ and $c \in (0,1)$ is the Hölder exponent of $f$.

Let $G(m,n)$ be the $Y$-block given by the concatenation

$$C[A_0, \langle]\ \langle^{K_i} \frown \mathbf{L}^{2N+m} \frown \mathbf{R}_i^{2N+n} \frown \rangle^{K_i}\ C[\rangle, A_1],$$

where $N, m, n$ are nonnegative integers. We first choose a suitable value for $N$. Consider

$$z(m,n) = \cdots \rangle\,\rangle \frown .\, A_0\ G(m,n)\ A_1 \frown \langle\,\langle \cdots$$

We think of the blocks $\mathbf{L}^{2N+m}$ as concatenations $\mathbf{L}^N\ \mathbf{L}^m\ \mathbf{L}^N$ and similarly with the $\mathbf{R}_i$ blocks. Writing it this way makes it clear, using (6), that the length in the suspension of the first $\mathbf{L}^N$ block is exponentially close (in $N$) to the corresponding block in $z(0,0)$ uniformly in $m$ (and $n$). Similarly with the second $\mathbf{L}^N$ block. Additionally, the length in the suspension of the $\mathbf{L}^m$ block is exponentially close in $N$ to $m \cdot \operatorname{PLen}(\overline{\mathbf{L}}, f)$, again uniformly in $m$ and $n$. Similar statements hold for the $\mathbf{R}_i$ blocks.

Similarly, the differences of the lengths in the suspension of the other corresponding blocks appearing in the concatenations forming $z(m,n)$ and $z(0,0)$ approach 0 exponentially in $N$, uniformly in $m$ and $n$.

Hence if the length in the suspension of the block $A_0\, G(0,0)$ in $(z(0,0),0)$ is $b(N)$ (the $N$ dependence appears in $G(0,0)$), the arguments above show that $N$ can be chosen so that the length in the suspension of the block $A_0\, G(m,n)$ in $z(m,n)$ is within $2^{-i}$ of $b(N) + m\operatorname{PLen}(\overline{\mathbf{L}}, f) + n\operatorname{PLen}(\overline{\mathbf{R}_i}, f)$ uniformly in $m$ and $n$. Let

$$\begin{aligned} b' \;=\;& \operatorname{MLen}(A_0, f) + \operatorname{MLen}\Big(C[A_0, \langle]\ \langle^{K_i} \frown \mathbf{L}^{2N}, f\Big) + \\ & \operatorname{MLen}\Big(C[\mathbf{L}, \mathbf{R}_i]\ \mathbf{R}_i^{2N}, f\Big) + \operatorname{MLen}\Big(C[\mathbf{R}_i, \rangle] \cdot \rangle^{K_i} \cdot C[\rangle, A_1], f\Big). \end{aligned}$$

Note that $b(N) \le b'$. By part (B) $\gcd\big(\operatorname{PLen}(\overline{\mathbf{L}}, f), \operatorname{PLen}(\overline{\mathbf{R}_i}, f)\big) \le 2^{-i}$, thus conditions (Fi) and (Fii) follow by choosing $L_i \ge b' + M_i - \operatorname{MLen}(A_0, f)$, where $M_i$ is such that for all $M > M_i$, there is a nonnegative integer combination of $\operatorname{PLen}(\overline{\mathbf{L}}, f)$ and $\operatorname{PLen}(\overline{\mathbf{R}_i}, f)$ that approximates $M$ to within an error of $2 \cdot 2^{-i}$.

To see that condition (Fiii) holds compare $z$ with $w'$. Observe that by (6) and the definition of $D$, the difference between the length in the suspension of $A_0 D_0$ in $z$ and $w'$ is bounded by a function of $i$ that goes to 0 as $K_i \to \infty$. Thus we may choose $K_i$ sufficiently large so that $A_1$ begins within $3 \cdot 2^{-i}$ of $L$ in $(w', 0)$ for all $w \in Y$. □

*Proof of Lemma 8* (H). Rather than deducing the required result from the basic definitions of topological entropy, we will use known results that are much more powerful.

Recall that $\mathbf{F}$ was defined in the proof of Lemma 8 (C) and (D) so that $h' := h_{\text{top}}(\text{susp}(Y \setminus \mathbf{F}), T^1) = h_{\text{top}}(X, \theta) + 2\delta$, for some $\delta > 0$. Let $P$ be the set of all periodic points of $Y \setminus \mathbf{F}$. By [3, Theorem 4.1] in conjunction with the isomorphism theorem of Bowen [4, Theorem 2], we have that for all $u > 0$ sufficiently large,

$$\kappa_1 \frac{e^{h'u}}{u} \leq \# \{y \in P : \text{PLen}(y, f) \leq u\} \leq \kappa_2 \frac{e^{h'u}}{u}, \tag{7}$$

for some constants $\kappa_1, \kappa_2 > 0$.

For the shift space $(X, \theta)$, we have that for $i$ sufficiently large

$$\#B^i(X) < e^{i(h'-\delta)}. \tag{8}$$

Note that by (6), for some constant $\kappa > 0$, the difference between the length in the suspension of a block and its maximum length in the suspension is at most $\kappa$. Thus by (7) and (8), we have that there exists $\alpha \in (0, 1)$, such that for all $i$ sufficiently large

$$\#\big(B^i(X) \cup B^{i+1}(X)\big) \leq \#B_{(1-\alpha)i}(Y \setminus \mathbf{F}, f).$$

Condition (Hii) is easily satisfied by taking $i$ sufficiently large, since the maximum lengths in the suspension of all connecting blocks are bounded above. Conditions (Hi) and (Hiii) can be verified with the definitions of $\mathbf{F}$, $\langle$, and $\rangle$. □

For a set $A$, let $\text{pow}(A)$ denote its power set.

**Lemma 9** (Markers upstairs)**.** *Let $X$ be a non-trivial subshift with an invariant measure $\mu$. Let $(\ell_i)_{i=1}^{\infty}$ be a sequence of positive integers such that $\ell_{i+1}/\ell_i \geq 4$. There exists a measurable function $M : \mathbb{Z}^+ \times X \to \text{pow}(\mathbb{Z})$ and a set $X' \subset X$ of full measure with the following properties.*

*(a) For all $i \in \mathbb{Z}^+$ and for all $x \in X'$, we have $M(i, \theta(x)) = M(i, x) - 1$.*
*(b) For all $0 < i < j$ and all $x \in X'$, we have $M(j, x) \subset M(i, x)$.*
*(c) For all $x \in X'$, the distance between two successive elements of $M(1, x)$ is $\ell_1$ or $\ell_1 + 1$.*
*(d) For all $i \geq 2$ and $x \in X'$, if $u < v$ are two consecutive elements of $M(i, x)$, then $\ell_i \leq v - u \leq 2\ell_i$ and the number of elements of $M(i - 1, x)$ in $[u, v]$ is between $\ell_i/(2\ell_{i-1})$ and $2\ell_i/\ell_{i-1}$.*

*Proof.* By Lemma 6, there exist two measurable sets $Q_1, Q_2$, such that $\{\theta^j(Q_1)\}_{j=0}^{\ell_1-1} \cup \{\theta^j(Q_2)\}_{j=0}^{\ell_1}$ give a partition of $X$ (modulo a null set). For each $x \in X$, we let

$$M(1,x) = \{u \in \mathbb{Z} : \theta^u(x) \in Q_1 \cup Q_2\}\,.$$

Clearly, conditions (a) (for $i = 1$) and (c) are satisfied.

Assume that $M(i,x)$ has been defined so that conditions (a), (b), and (d) hold for all $j \leq i$. It remains to define $M(i+1,x)$. Again, by Alpern's theorem, there exist two measurable sets $P_1, P_2$, such that $\{\theta^j(P_1)\}_{j=0}^{\ell_{i+1}+2\ell_i-1} \cup \{\theta^j(P_2)\}_{j=0}^{\ell_{i+1}+2\ell_i}$ give a partition of $X$. Consider the set

$$M(i+1,x)' = \{u \in \mathbb{Z} : \theta^u(x) \in P_1 \cup P_2\}\,.$$

For each $u' \in M(i+1,x)'$, let $\kappa(u') = \inf\{u \in M(i,x) : u \geq u'\}$. Set $M(i+1,x) = \{\kappa(u') : u' \in M(i+1,x)'\}$. □

The following definitions will be important in defining the sequence of maps in Proposition 7. Let $M$ be the function from Lemma 9. For each $i \in \mathbb{Z}^+$ and each $x \in X$, let

$$N_i^-(x) := \sup\{z \in M(i,x) : z \leq 0\} \qquad (9)$$

$$\text{and } N_i^+(x) := \inf\{z \in M(i,x) : z > 0\}\,. \qquad (10)$$

**Corollary 10.** *Let $(a_i)_{i=1}^\infty$ be a sequence of positive integers. If*

$$\ell_{i+1} \geq 2^{i+2} a_i \ell_i, \qquad (11)$$

*then the function $M$ from Lemma 9 has the additional property that*

$$\mu(\{x \in X : N_{i+1}^+(x) - N_i^+(x) \leq 2a_i\ell_i\}) \leq 2^{-i}; \text{ and}$$
$$\mu(\{x \in X : N_i^-(x) - N_{i+1}^-(x) \leq 2a_i\ell_i\}) \leq 2^{-i}.$$

*Hence for $\mu$-almost every $x \in X'$, $N_i^+(x) \to \infty$ and $N_i^-(x) \to -\infty$ as $i \to \infty$.*

*Proof.* We estimate $\mu(\{x \in X\colon N_{i+1}^+(x) - N_i^+(x) \leq 2a_i\ell_i\})$, the other estimate being identical. From Lemma 9, $N_i^+(x)$ is less than $2\ell_i$ for all $x$. Hence $N_{i+1}^+(x) - N_i^+(x) \leq 2a_i\ell_i$ implies $N_{i+1}^+(x) < 2(a_i+1)\ell_i \leq 4a_i\ell_i$.

Next, notice since consecutive $(i+1)$-markers are separated by at least $\ell_{i+1}$, one has for each $0 \leq k < \ell_{i+1}$, $N_{i+1}^+(x) = \ell_{i+1}$ if and only if $N_{i+1}^+(\theta^k(x)) = \ell_{i+1} - k$. Hence the sets $(N_{i+1}^+)^{-1}(k)$ have equal measure for $k$ in the range $1 \leq k \leq \ell_{i+1}$ so that $\mu(\{x\colon N_{i+1}^+(x) = j\}) \leq 1/\ell_{i+1}$ for all $j \leq \ell_{i+1}$. We then see that $\mu(\{x\colon N_{i+1}^+(x) \leq 4a_i\ell_i\}) \leq 4a_i\ell_i/\ell_{i+1} < 2^{-i}$ as required. □

## 4. Definition of Injection

*Proof of Proposition 7.* Suppose that $(X, \theta)$ is a non-trival ergodic subshift with invariant measure $\mu$, and $\mathrm{susp}(Y, f)$ is a topologically mixing suspension flow whose topological entropy exceeds the topological entropy of $(X, \theta)$ as in the statement of the proposition. Let $W$ be a $Y$-block. Let $\langle, \rangle, \mathbf{L}, \mathbf{R}_1, \mathbf{R}_2 \ldots,$ $(K_n, k_n, L_n)_{n\in\mathbb{Z}^+}$, and $\alpha \in (0, 1)$ be given by Lemma 8 ($K_n$ being the number of $\langle$ and $\rangle$ used to delimit the filler blocks between level $n$ blocks, $L_n$ being a length in the suspension such that any gap of size exceeding $L_n$ can be filled up to within $2^{-n}$ accuracy by level $n$ filler blocks, and $\alpha$ being the compression that can be achieved in encoding the first level information blocks). The block $W$ appears in $\langle$ and $\rangle$.

Note that by (6), for some constant $\kappa > 0$, the difference between the length in the suspension of a block and its maximum length in the suspensions is at most $\kappa$. Let $(a_i)_{i\in\mathbb{Z}^+}$ be a sequence of positive integers such that

$$a_i(L_i - \kappa) > 2L_{i+1}. \tag{12}$$

The $a_i$ represent the number of level $i$ blocks forming the crumple zone when building the level $i+1$ blocks. Choose $\ell_1 > 0$ large enough so that by Lemma 8 (H), there exists an injection $\phi : B^{\ell_1+1}(X) \cup B^{\ell_1}(X) \to B_{(1-\alpha)\ell_1}(Y \setminus \mathbf{F}, f)$ such that $\rangle\, \phi(B) \,\langle \in B(Y)$ for all $B \in B^{\ell_1+1}(X) \cup B^{\ell_1}(X)$. We can additionally require that $\ell_1$ is large enough that

$$\alpha\ell_1 > 2L_1. \tag{13}$$

We define the subsequent $\ell_i$ by setting

$$\ell_{i+1} = 2^{i+2}a_i\ell_i. \tag{14}$$

Let $M$ and $X' \subset X$ be given be Lemma 9, and $N_i^+$ and $N_i^-$ be given by (9) and (10), respectively. By Corollary 10, $N_i^+(x) \to \infty$ and $N_i^-(x) \to -\infty$ for $\mu$-almost all $x \in X'$.

We will recursively define sequences of maps $\phi_i : X' \to B(Y)$. Each $\phi_i$ will be extended to a map $\psi_i \colon X' \to Y$ by the relation

$$\psi_i(x) = \cdots \rangle \,\rangle\mathbf{.}\phi_i(x) \,\langle\, \langle \cdots . \tag{15}$$

We set

$$\Psi_i(x) := T^{-N_i^-(x)}(\psi_i(x), 0) \tag{16}$$

for all $x \in X'$. Given $(y, s) \in \mathrm{susp}(Y, f)$, let $\pi_1(y, s) = y$ and $\pi_2(y, s) = s$. We shall inductively construct maps $\phi_i, \psi_i, \Psi_i$ satisfying the following properties (previously illustrated in Figure 1).

(I) For all $x \in X'$; and for all $i \geq 1$, we have $\rangle\, \phi_i(x) \,\langle \in B(Y)$.

(II) For all $x \in X'$ and for all $i \geq 1$, we have

$$N_i^+(x) - N_i^-(x) > 2L_i + \mathrm{MLen}(\phi_i(x), f).$$

(III) For all $x \in X'$, $i \geq 1$, and all $k \in [N_i^-(x), N_i^+(x))$, we have

$$\Psi_i(\theta^k x) = T^k \Psi_i(x).$$

(IV) For all $x \in X'$ and for all $i \geq 2$, if $|N_i^+(x) - N_{i-1}^+(x)| > 2a_i\ell_i$, then there exists $|\varepsilon| \leq 3 \cdot 2^{-i}$ such that

$$\pi_1\big(T^{N_{i-1}^-(x)+\varepsilon}\Psi_i(x)\big)_0^{n-1} = \phi_{i-1}(x) = \pi_1\big(T^{N_{i-1}^-(x)}\Psi_{i-1}(x)\big)_0^{n-1}; \text{ and}$$
$$\pi_2\big(T^{N_{i-1}^-(x)+\varepsilon}\Psi_i(x)\big) = 0 = \pi_2\big(T^{N_{i-1}^-(x)}\Psi_{i-1}(x)\big),$$

where $n$ is the size of $\phi_{i-1}(x)$.

Let $x \in X'$. Set $\phi = \phi_1$. By Lemma 8 (H) (Hii), we may set

$$\psi_1(x) := \cdots \rangle \rangle \,.\phi_1(x_{N_1^-(x)} \cdots x_{N_1^+(x)-1}) \langle \langle \cdots .$$

Note that by Lemma 9 (c), we have that $N_1^+(x) - N_1^-(x) \in \{\ell, \ell+1\}$. Clearly, by (13) and the definition of $\phi$ conditions (I), (II) and (III) are satisfied.

Suppose $\phi_i$ and $\psi_i$ have been defined and satisfy conditions (I), (II), (III), and (IV). Consider the set

$$\left\{v \in M(i,x) : v \in [N_{i+1}^-(x), N_{i+1}^+(x))\right\}.$$

Let $N_{i+1}^-(x) = n^1 < n^2 < \cdots < n^{j-1} < n^j < N_{i+1}^+(x)$ be an enumeration of the set. Let $j' = j - a_i$ (this is positive because the number of level $i$ blocks in a level $i+1$ block is at least $\ell_{i+1}/(2\ell_i) > 2^{i+1}a_i$ by Lemma 9 (d) and by (14)). The point $\psi_{i+1}(x)$ will be of the form

$$\cdots \rangle \rangle \,.B_1F_1B_2F_2 \cdots F_{j-1}B_j \langle \langle$$

where $B_k = \phi_i(\theta^{n^k}x)$ for each $1 \leq k \leq j$ and the $F_k$ are given by $F(\ell_k, r_k)$, where

$$F(\ell, r) := \langle^{K_i} \frown \mathbf{L}^\ell \frown \mathbf{R}_i^r \frown \rangle^{K_i} \,. \tag{17}$$

We call such an $F$ a ***level* i *filler block***.

By repeated applications of Lemma 8 we can choose the $\ell_k$ and $r_k$ such that for $1 \leq k \leq j'$, $B_k$ begins within $3 \cdot 2^{-i}$ of $n^k$ in $\Psi_{i+1}(x) = T^{-n^1}\psi_{i+1}(x)$. For $k \in [j', j)$ we simply set $F_k = F(1,1)$ to create the crumple zone mentioned earlier in Figure 1.

We set $\phi_{i+1}(x) = B_1F_1B_2 \cdots F_{j-1}B_j$, and note that property (I) is trivially satisfied. We call the $\phi_{i+1}(x)$ ***level* i + 1 *blocks***.

Notice that the crumple zone consists of $a_i$ level $i$ blocks. Each of these consists of at most $2\ell_i$ symbols by Lemma 9 (d), so that $|N_i^+(x) - N_{i-1}^+(x)| > 2a_i\ell_i$ guarantees that $x$ lies in a level $i$ block outside the crumple zone. The construction above then ensures that (IV) holds for $i+1$.

By construction, we have that $\mathrm{MLen}(F_k, f) \leq L_i + \kappa$ for $k \in [j', j)$, so using (II) at level $i$ for these $k$'s we have

$$n^{k+1} - n^k - \mathrm{MLen}(F_k, f) \geq L_i - \kappa.$$

Since there are $a_i$ such $k$'s, we obtain

$$N^+_{i+1}(x) - N^-_{i+1}(x) - \mathrm{MLen}(\phi_{i+1}(x), f) \geq a_i(L_i - \kappa)$$

Hence by (12), condition (II) is satisfied at level $i+1$. Note that property (III) follows immediately from the definition of $\Psi_{i+1}$. Hence the inductive step is complete.

Conclusion (a) (of Proposition 7) follows immediately from property (III) and Corollary 10. Conclusion (b) follows from property (IV), Corollary 10, and the Borel-Cantelli Lemma.

It remains to show that properties (c), (d), and (e) (of Proposition 7) hold. Suppose that

$$\Psi(x) := \lim_{i\to\infty} \Psi_i(x) = (\mathbf{y}(x), \mathbf{s}(x)) \in \mathrm{susp}(Y, f),$$

for all $x \in X'' \subset X'$, where $X''$ is also a set of full measure. Let $x \in X''$. Let $(n_i)_{i\in\mathbb{Z}}$ be an enumeration of the set $M(1, x)$, where $n_0 = N^-_1(x)$. For each $i \in \mathbb{Z}$, let

$$B_i := \varphi(x_{n_i} \cdots x_{n_{i+1}-1}).$$

Recall that in the construction, level 1 blocks (images of $\phi$) are combined into longer level 2 blocks by placing level 1 filler between the pairs of blocks while the leftmost level 1 block in a level 2 block does not (yet) have filler next to it and similarly for the rightmost level 1 block in a level 2 block. At the next stage level 2 blocks are interspersed with level 2 filler, again leaving the end blocks bare. Hence, as long as $N^+_i(x) \to \infty$ and $N^-_i(x) \to \infty$, then the level of the filler on the left of $\phi(x)$ is $\min\{k\colon N^-_{k+1}(x) \neq N^-_1(x)\}$ with a similar expression for the level of filler on the right. Hence for each $x \in X''$, $\mathbf{y}(x)$, the first coordinate of $\Psi(x)$, is of the form

$$\cdots B_{-1}F_{-1}B_0F_0B_1F_1\cdots, \tag{18}$$

where the $B_i$ are of the form $\phi(C)$ for $C \in B^{\ell_i}(X) \cup B^{\ell_i+1}(X)$ and the $F_i$ are $Y$-blocks of the form $\langle^{K_n} C \rangle^{K_n}$ for some $n$ and some $Y$-block $C$; furthermore, by Lemma 8 (D) and (G), we have that the block $C$ does not contain a $\langle^{K_j}$ or $\rangle^{K_j}$, with $j \geq n$.

Note that the blocks $\rangle^n$ and $\langle^n$ appear in $\mathbf{y}(x)$ for all $n \in \mathbb{Z}^+$. Hence by Lemma 8 (D), property (c) holds. Since $W$ appears in $\langle$, property (e) holds.

Let $y \in Y$ be of the form (18). If $y_0$ belongs to one of the $F$ blocks, then one of the following happens:

(1) there exist $m \leq 0 < n$ such that $y_m^{n-1}$ is $\langle$ or $\rangle$; Or
(2) there exists an $i \geq 1$ and an $n > 0$ such that $\rangle^{K_i}$ occurs in $y_0^{n-1}$ but $\langle^{K_i}$ does not.

On the other hand, if $y_0$ belongs to one of the $B$ blocks, then (2) is ruled out by the properties of the filler blocks. (1) is ruled out by properties of the level one blocks in Lemma 8 (Hi) and (Hiii).

For each symbol in $y$, one can therefore decide whether it belongs to one of the $B$ blocks or one of the $F$ blocks. Hence given the first coordinate of $\Psi(x)$, (the element of $Y$) one can recover the sequence of $B_i$ blocks and hence applying $\phi_1^{-1}$ to each of them, one can recover $x$ up to translation, completing the proof of property (d).

□

## 5. Proof of Theorem 1

Bowen and Ratner constructed symbolic dynamics for geodesic flows, proving the following theorem.

**Theorem 11** (Bowen, R. [5] and Ratner, M. [31])**.** *Let $(W^t)_{t\in\mathbb{R}}$ be geodesic flow on a compact surface of variable negative curvature $\mathcal{M}$, with unit tangent bundle* $\mathrm{UT}(\mathcal{M})$. *There exists a suspension flow* $(\mathrm{susp}(Y, f), (T^t)_{t\in\mathbb{R}})$ *over a subshift of finite type $Y$ under a Hölder continuous function $f : Y \to (0, \infty)$, and a finite-to-one continuous surjection from $\pi : \mathrm{susp}(Y, f) \to \mathrm{UT}(\mathcal{M})$ such that $\pi \circ W^t = T^t \circ \pi$ for all $t \in \mathbb{R}$.*

Furthermore, the subshift of finite type $Y$ can be chosen to be irreducible [7, Lemma 2.1], [6, Section 3], and $\pi$ is one-to-one on a set of full measure for any measure that is ergodic and fully supported [29, Theorem III.8].

Note that if map $\pi$ in Theorem 11 were one-to-one for all ergodic measures, then Theorem 1 would follow immediately from Theorem 2. We will not be able to use Theorem 11 directly because the map we define for Theorem 3 has an image which is in general not fully supported.

One approach is to modify our proof of Theorem 3 so that the embedding has an image that is fully supported; this may be accomplished by enumerating the countable number of $Y$-blocks, and inserting them into the filler blocks so that each $Y$ block appears in some level $i$ filler block, for some $i \geq 1$. Carrying out this more technical construction has the advantage that we obtain full universality, rather than just universality. Instead, we will use the following corollary, which is a consequence of the Proof of Theorem 11 given by Bowen [5].

**Corollary 12** (Corollary of Bowen's proof of Theorem 11)**.** *There exists $W \in B^m(Y)$ for some $m$ such that $\pi$ is one-to-one on the set of all $(y, s) \in \mathrm{susp}(Y, f)$ such that for infinitely many positive $n \geq 0$ and infinitely many $n < 0$, we have $(\theta^n y)_0^{m-1} = W$.*

Corollary 12 will also be useful in providing negative answer to the second question of about geodesic flows due to Ledrappier and Federico and Jana Rodriguez Hertz.

*Proof of Corollary 12.* Let $\mathcal{P}$ be the Markov partition of $UT(\mathcal{M})$ provided by [5, Section 2 and Lemma 7.5]; the elements of $\mathcal{P}$ are the closures of their interiors, and that their interiors are disjoint. Let $B \in \mathcal{P}$. By expansiveness of geodesic flow, there exists a cylinder set $W$ in $Y$ and an interval $J$ such that $\pi(W \times J) \subset \mathrm{int}\, B$. The map $\pi$ is one-to-one on the preimage of the set of points whose orbits never intersect the stable or unstable boundaries of the Markov partition. Further, the union of the stable boundaries of the elements of the partition is forwards-invariant, while the union of the unstable boundaries of the elements of the partition is backwards-invariant.

If $W$ appears infinitely often in the future of $y$, then $\pi(y, s)$ cannot belong to the stable boundary. If $W$ appears infinitely often in the past of $y$, then $\pi(y, s)$ cannot belong to the unstable boundary. By the above observations about invariance, the orbit of $\pi(y, s)$ never hits the stable or unstable boundaries of the partition and so $\pi^{-1}(\pi(y, s)) = \{(y, s)\}$ as required. □

*Proof of Theorem 1.* As in the proof of Theorem 2, by the Jewett-Krieger theorem [20], we may make the simplifying assumption that the measure-preserving system that is to be embedded is a uniquely ergodic subshift $(X, \mu, \theta)$. Let $(\mathrm{susp}(Y, f), (T^t)_{t\in\mathbb{R}})$, $\pi$, and $W$ be given by Theorem 11 and Corollary 12. Now repeating verbatim the proof of Theorem 3, but using the specific word $W$ gives and embedding of $(X, \theta, \mu)$ into $(\mathrm{susp}(Y, f), T^1)$, where $W$ appears in first coordinate of the image points infinitely many times in the future and the past. Corollary 12 then implies that $\pi$ is one-to-one on the image of $\Psi$. The composition yields the result. □

*Proof of Corollary 4.* Let $(\Omega', \mathcal{D}', \nu')$ be a non-atomic probability space endowed with an ergodic measure-preserving automorphism $S'$. Consider the automorphism $S$ defined on $\Omega := \Omega' \times \{0, 1\}$ given by $S(x, 0) := (x, 1)$ and $S(x, 1) := (S'(x), 0)$ for all $x \in \Omega'$, and the measure $p$ on $\{0,1\}$ such that $p(0) = p(1) = 1/2$. It is easy to verify that $S$ is ergodic and preserves the measure $\nu := \nu' \times p$. Furthermore, by [11, Lemma 8.7], $S$ does not have a square-root; that is, there does not exist a subset of $\Omega \times \{0,1\}$ of full measure for which there is measure-preserving automorphism $U$ such that $U \circ U = S$ on the subset. Hence the result follows immediately from Theorem 1 and choosing $S'$ and thus $S$ with measure-theoretic entropy sufficiently small. □

## 6. Appendix

The arguments given above also allow us to resolve the second question about geodesic flows due to Ledrappier and Federico and Jana Rodriguez Hertz, who asked if a minimal subset for the $\mathbb{R}$-action of geodesic flow (without periodic points of rational period) is necessarily minimal for the $\mathbb{Z}$-action of the time-one map of the geodesic flow. We give a counterexample to this using Corollary 12 above.

**Proposition 13.** *Let $(G^t)_{t\in\mathbb{R}}$ be the geodesic flow on a compact surface of variable negative curvature $\mathcal{M}$, with unit tangent bundle $\mathrm{UT}(\mathcal{M})$. Then there exists a minimal subset $C$ for the $\mathbb{R}$-action on $\mathrm{UT}(\mathcal{M})$ that is not minimal for the $\mathbb{Z}$-action.*

*Proof.* As before, let $\mathrm{susp}(Y, f)$ be a suspension over an irreducible shift of finite type and $\pi$ be a continuous factor map from the suspension flow to the geodesic flow. By Corollary 12, let $W$ be a $Y$-block such that on the set of points for which $W$ appears infinitely often in the past and the future, $\pi$ is one-to-one.

If $Y$ contains a periodic point $\overline{A}$ with rational period in the suspension, then let $x = \pi(\overline{A}, 0)$. Clearly the $\mathbb{R}$-orbit of $x$ is closed and the $\mathbb{Z}$-orbit is a discrete subset.

Suppose then that $\mathrm{susp}(Y, f)$ has no periodic points with rational period in the suspension. Let $B_0$ and $B_1$ be two $Y$-blocks in which $W$ appears, such that $B_i\,B_j$ is a $Y$-block for all pairs $i, j \in \{0,1\}$. Let them agree on sufficiently many symbols at each end that switching any number of $B_0$'s for $B_1$'s or vice versa in a point $y \in Y$ does not change the length in the suspension of any contiguous block (disjoint from the blocks being switched) by more than $\frac{1}{10}$ (see (6)). We may also assume that $\overline{B_0}$ does not contain any $B_1$'s and $\overline{B_1}$ does not contain any $B_0$'s. Finally we will assume that the $B_i$ are not powers of smaller words. See Lemma 8 (D) for a similar construction.

We build a point $y$ in $Y$ in the following way: Let

$$y^{(n)} = \cdots B_0\,B_0\boldsymbol{.}\,B_{i_1} \cdots B_{i_n}\ B_0\,B_0 \cdots$$

where the $(i_n)$ are defined recursively by

$$i_n = \begin{cases} 1 & \text{if the block } B_{i_{n-1}} \text{ ends within } \frac{1}{10} \text{ of } \mathbb{Z} \text{ in } y^{(n-1)}; \\ 0 & \text{otherwise.} \end{cases}$$

The limit point is then $y = \cdots B_0\, B_0 . B_{i_1}\, B_{i_2} \cdots$. Notice that in $y$ the $B_0$ and $B_1$ both appear with bounded gaps (because in any sufficiently long block of $B_0$'s, say, their lengths in the suspension become arbitrarily close to $\mathrm{PLen}(\overline{B_0}, f)$, the fractional parts of whose multiples are dense in $[0,1)$). Also all $B_1$'s start at points of $\mathbb{Z} + [-\frac{1}{5}, \frac{1}{5}]$ (when subsequent blocks are altered the starting points may move by up to $\frac{1}{10}$). Let $\Omega$ denote the $\omega$-limit set of $y$ under the $\mathbb{R}$-action. By the above observation, all points of $\Omega$ have $B_0$'s and $B_1$'s appearing with bounded gaps. Further, for each $(z,s) \in \Omega$, there is a $\beta \in [0,1)$ such that all $B_1$'s begin at points of $\mathbb{Z} + \beta + [-\frac{1}{5}, \frac{1}{5}]$. Let $C$ be a minimal subset of $\Omega$ (under the $\mathbb{R}$-action) and let $(z,s) \in C$. Suppose that $B_1$'s start at points of $\mathbb{Z} + \beta + [-\frac{1}{5}, \frac{1}{5}]$. Then the same is true of the $\mathbb{Z}$-orbit closure of $(z,s)$, $\Omega' \subseteq C$. So $\Omega'$ is a closed $\mathbb{Z}$-invariant subset of $C$. Since $T^{\frac{1}{2}}(z,s) \notin \Omega'$, we see that $\Omega'$ is a proper closed $\mathbb{Z}$-invariant subset of $C$, so that $C$ is not minimal for the $\mathbb{Z}$-action.

By Corollary 12, since $C$ consists of points with infinitely many $W$'s in the past and the future, the restriction of $\pi$ to $C$ is one-to-one. Since $C$ is compact, $\pi$ is a conjugacy between the real action on $C$ and the real action on $\pi(C)$. Hence (lack of) minimality is preserved and $\pi(C)$ is minimal for the $\mathbb{R}$-action on $\mathrm{UT}(\mathcal{M})$ but not for the $\mathbb{Z}$ action. □

In the following example, we construct a filling function over the full shift taking only rational values. Further it has the property that $\mathrm{susp}(X \setminus F, f)$ fails to be topologically weak mixing for any word $F$.

*Example* 1. Let $X$ be the full two shift on $\{0,1\}$. Let $N_k = 3^k$, $E_k = \{x : x_0^{N_k - 1} \text{ contains all } k\text{-blocks}\}$ and $a_k = 1 + 2^{-N_{k+1}}$. Let $K(x) = \min\{k : x \notin E_k\}$ and define $f(x) = a_{K(x)}$.

This is a Lipschitz function taking only rational values. Now fix any word $w$, of size $k$, say. All points of $X \setminus w$ take $f$ values that are a multiple of $2^{-N_{k+1}}$ so that the sums of $f$ over periodic points that don't contain any $w$'s are multiples of $2^{-N_{k+1}}$.

However the greatest common divisor of the full set of periods is 0. To see this, let $w_1$ and $w_2$ be two distinct words of size $k$. Assume further that $0^k w_1$ does not contain a $w_2$ and $w_2 0^k$ does not contain a $w_1$. Take a block of size $N_k$ containing $w_1$ as its first $k$ symbols and $w_2$ as its last $k$ symbols with no $w_1$'s or $w_2$'s between. Further ensure that the block contains all $k$ words. Extend this to a block of size $2N_k$ by adding $N_k$ 0's and let $x$ be the periodic orbit obtained by concatenating this word. Summing $f$ over the period, one obtains exactly one value of $1 + 2^{-N_{k+2}}$ while all the other values are multiples of $2^{-N_{k+1}}$ so that the period in the suspension is an odd multiple of $2^{-N_{k+2}}$. Since the point $\overline{0}$ has a period in the suspension of $1 + 2^{-9}$, the greatest common divisor of these two periods is a factor of $513/2^{N_{k+2}}$. Since this holds for all $k$, the example is complete.

## Acknowledgements

We thank Jean-Paul Thouvenot for generous suggestions and Michael Hochman, Mark Pollicott, and Richard Sharp for helpful comments. We also thank the referees for their careful reading of this paper and useful suggestions.

(A. Quas and T. Soo) Department of Mathematics and Statistics, University of Victoria, PO BOX 3060 STN CSC, Victoria, BC V8W 3R4, Canada
*E-mail address*: aquas at uvic.ca
*URL*: http://www.math.uvic.ca/~aquas/

*E-mail address*: tsoo at uvic.ca
*URL*: http://www.math.uvic.ca/~tsoo/